\documentclass{article}
\usepackage{amssymb,amsmath}
\usepackage{amsfonts}
\usepackage{epsfig}

\begin{document}

\begin{center}
{\Large Holomorphic functions in generalized\\[0pt]
\vspace{3mm} Cayley-Dickson algebras}
\end{center}

\vspace{5mm}

\begin{center}
{\large Cristina Flaut and Vitalii Shpakivskyi}
\end{center}

\vspace{5mm}

\textbf{Abstract.} {\small In this paper we investigated some properties of
holomorphic functions (belonging to the kernel of the Dirac operator)
defined on domains of the real Cayley-Dickson algebras. For this purpose, we
study first some properties of these algebras, especially multiplication
tables for certain elements of the basis. Using these properties, we
provided an algorithm for constructing examples of the class of functions
under consideration.\medskip }%
\begin{equation*}
\end{equation*}%
{\small \ }

\textbf{Keywords:} Cayley-Dickson and generalized Cayley-Dickson algebras;
Dirac operator; holomorphic functions. \vspace{5mm}

\textbf{AMS Subject Classification:} 30G35, 17A20. \vspace{3mm}%
\begin{equation*}
\end{equation*}

\textbf{0. Introduction}\vspace{5mm}

The theory of quaternionic differentiable functions has many applications in
different areas of mathematics, physics and in other applied sciences (see,
for example, \cite{Gu-Shp-97}, \cite{Kr-Sh-96}). This theory has its origins
in the paper \cite{Mojsil-Teodoresku} in which the authors proposed, for the
first time, an analogue of the Cauchy-Riemann conditions in
three-dimensional case. For the four-dimensional case, an analogue of these
conditions was considered in the paper \cite{Fu-35} and, as a next step of
this generalization, the differentiable functions in the octonionic algebra
was considered in the papers \cite{Lim-Shon-13}, \cite{Xi-Zh-Li-05}.

Generalization of the Cauchy-Riemann conditions in all algebras obtained by
the Cayley-Dickson process (called Cayley-Dickson algebras) was done in the
paper \cite{Ludkovsky-07}, where differentiable functions of variables
belonging to Cayley-Dickson algebras were defined. For such functions, was
established analogues results with the main results of complex analysis,
results which can be successfully used in the further studies of special
functions of variables with values in Cayley--Dickson algebras.

Comparing with \cite{Ludkovsky-07}, in the present paper, we investigate
another class of differentiable functions (using the Dirac operator) in
Cayley-Dickson algebras and, more important, we provide an example of this
kind of functions and an algorithm to find such as examples. Since these
functions are rather complicated objects, it is quite important to have a
way to generate examples.

The paper is organized in two sections. In the\ first section, we briefly
presented some properties of algebras obtained by the Cayley-Dickson process
and the algorithm described by \ J. W. Bales regarding to an easy way to
multiply the elements from a basis in such algebras (by using \textit{%
exclusive or }operation and a \textit{twist map). }In the second section, by
description the multiplication tables for certain elements of the basis
(Propositions 2.2 and 2.3), we obtained the main result of this work: an
example of a left\textit{\ }hyperholomorphic function \ in generalized \
Cayley-Dickson algebras (Theorem 2.12). Moreover, in the Thorem 2.10 we
proved that \ for studying left $A_{t}$-holomorphic functions in generalized
Cayley-Dickson algebras $A_{t}=\left( \frac{\gamma _{1},...,\gamma _{t}}{%
\mathbb{R}}\right) $ it is suffices to consider left $A_{t}$-holomorphic
functions only in the algebras $\left( \frac{sign(\gamma
_{1}),...,sign(\gamma _{t})}{\mathbb{R}}\right) ,$ that means this study
depends only by the sign of the real numbers $\gamma _{1},\gamma
_{2},...,\gamma _{t}.\medskip $

\begin{equation*}
\end{equation*}

\textbf{1. Preliminaries} \vspace{3mm}

Let $K$ be a commutative field with $charK\neq 2$ and $A$ be an algebra over
the field $K.~$A unitary algebra $A\neq K$ such that we have $x^{2}+\alpha
_{x}x+\beta _{x}=0,$ for each $x\in A,$ with $\alpha _{x},\beta _{x}\in K,$
is called a \textit{quadratic algebra}.

In the following, we briefly present the \textit{Cayley-Dickson process} and
the properties of the algebras obtained. For details about the
Cayley-Dickson process, the reader is referred to \cite{Sc-66} and \cite%
{Sc-54}.

Let $A$ be a finite dimensional unitary algebra over a field $\ K$ with a 
\textit{scalar} \textit{involution} 
\begin{equation*}
\,\,\,\overline{\phantom{xx}}:A\rightarrow A,\quad a\rightarrow \overline{a},
\end{equation*}%
$\,\,$ i.\thinspace e. a linear map satisfying the following relations:$%
\,\,\,\,\,$ 
\begin{equation*}
\overline{ab}=\overline{b}\overline{a},\quad \overline{\overline{a}}=a,
\end{equation*}%
and 
\begin{equation*}
a+\overline{a},a\overline{a}\in K\cdot 1\ \text{for all }a,b\in A.
\end{equation*}%
The element $\,\overline{a}$ is called the \textit{conjugate} of the element 
$a,$ the linear form$\,\,$ 
\begin{equation*}
\,\,t:A\rightarrow K,\quad t\left( a\right) =a+\overline{a}
\end{equation*}%
and the quadratic form 
\begin{equation*}
n:A\rightarrow K,\quad n\left( a\right) =a\overline{a}\ 
\end{equation*}%
are called the \textit{trace} and the \textit{norm \ }of \ the element $a.$
Hence an algebra $A$ with a scalar involution is quadratic. $\,$

Let $\gamma \in K$ \thinspace be a fixed non-zero element. We define the
following algebra multiplication on the vector space 
\begin{equation}
A\oplus A:\left( a_{1},a_{2}\right) \left( b_{1},b_{2}\right) :=\left(
a_{1}b_{1}+\gamma b_{2}\overline{a_{2}},\overline{a_{1}}b_{2}+b_{1}a_{2}%
\right) .  \label{1.1}
\end{equation}%
We obtain an algebra structure over $A\oplus A,$ denoted by $\left( A,\gamma
\right) $ and called the \textit{algebra obtained from } $A$ \textit{\ by
the Cayley-Dickson process} or simply \textit{generalized Cayley-Dickson
algebra}. We have $\dim \left( A,\gamma \right) =2\dim A$.

Let $x\in \left( A,\gamma \right) $, $x=\left( a_{1},a_{2}\right) $. The map 
\begin{equation*}
\,\,\,\overline{\phantom{x}}:\left( A,\gamma \right) \rightarrow \left(
A,\gamma \right) \,,\,\,x\rightarrow \bar{x}\,=\left( \overline{a}_{1},\text{%
-}a_{2}\right) ,
\end{equation*}%
is a scalar involution of the algebra $\left( A,\gamma \right) $, extending
the involution $\overline{\phantom{x}}\,\,\,$of the algebra $A.$

\thinspace If we take $A=K$ \thinspace and apply this process $t$ times, $%
t\geq 1,\,\,$we obtain an algebra over $K,\,\,$%
\begin{equation*}
A_{t}=\left( \frac{\gamma _{1},...,\gamma _{t}}{K}\right) .
\end{equation*}
By induction in this algebra, the set $\{e_{0}=1,e_{1},...,e_{n-1}
\},n=2^{t}, $ generates a basis with the properties: 
\begin{equation}  \label{1.2}
e_{i}^{2}=\gamma _{i}1,\,\,\gamma_{i}\in K,\gamma _{i}\neq 0,\,\,i=1,...,n-1
\end{equation}
and 
\begin{equation}  \label{1.3}
e_{i}e_{j}=-e_{j}e_{i}=\beta _{ij}e_{k},\,\,\beta _{ij}\in K,\,\,\,\beta
_{ij}\neq 0,\,\,\,i\neq j,\,\,\,i,j=\,\,1,...n-1,
\end{equation}
$\beta _{ij}$ and $e_{k}$ being uniquely determined by $e_{i}$ and $e_{j}.$

From \cite{Sc-54}, Lemma 4, it results that in an algebra $A_{t}$ with the
basis $B=\{e_{0}=1,e_{1},...,e_{n-1}\}$ satisfying relations (\ref{1.2}) and
(\ref{1.3}) we have:

\begin{equation}  \label{1.1}
e_{i}\left( e_{i}x\right) =\gamma _{i}^{2}x=(xe_{i})e_{i},
\end{equation}
for all $i\in \{1,2,...,n-1\}$ and for every $x\in A.$

The algebras $A_{t}$, in general, are neither commutative and nor
associative algebras, but are \textit{flexible} (i.\thinspace e. $%
x(yx)=(xy)x=xyx$, for all $x,y\in A_{t})$ quadratic and \textit{power
associative }(i.\thinspace e. the subalgebra $<x>$ of $A$, generated by any
element $x\in A$, is associative).

\vspace{2mm}

\textbf{Remark 1.1.} For $\gamma _{1}=...=\gamma _{t}=-1$ and $K=\mathbb{R},$
in \cite{Ba-09}, the author described \ how we can multiply the basis
vectors in the algebra $A_{t},\dim A_{t}=2^{t}=n$. He used the binary
decomposition for the subscript indices.

Let $\ e_{p},e_{q}$ be two vectors in the basis $B$ with $p,q$ representing
the binary decomposition for the indices of the vectors, that means $p,q$
are in $\mathbb{Z}_{2}^{n}.$ We have that $e_{p}e_{q}=\gamma _{n}\left(
p,q\right) e_{p\otimes q},$ where:

i) $p\otimes q$ are the sum of \ $p$ and $q$ in the group $\mathbb{Z}
_{2}^{n} $ or, more precisely, the "\textit{exclusive or}" for the binary
numbers $p$ and $q;$

ii) $\gamma _{n}$ is a function $\gamma _{n}:\mathbb{Z}_{2}^{n}\times 
\mathbb{Z}_{2}^{n}\rightarrow \{-1,1\}.$

The map $\gamma _{n}$ is called the \textit{twist map}.

The elements of the group $\mathbb{Z}_{2}^{n}$ \ can be considered as
integers from $0$ to $2^{n}-1$ with multiplication "\textit{exclusive or}"
of the binary representations. Obviously, this operation is equivalent with
the addition \ in $\mathbb{Z}_{2}^{n}.\medskip $

From now on, in whole the paper, we will consider $K=\mathbb{R}.$ Using the
same notations as in the Bales's paper, we consider the following matrices: 
\begin{equation}
A_{0}=A=\left( 
\begin{array}{cc}
1 & 1 \\ 
1 & -1%
\end{array}%
\right) ,\quad B=\left( 
\begin{array}{cc}
1 & -1 \\ 
1 & 1%
\end{array}%
\right) ,\quad C=\left( 
\begin{array}{cc}
1 & -1 \\ 
-1 & -1%
\end{array}%
\right) .  \label{1.6}
\end{equation}

In the same paper \cite{Ba-09}, the author find the properties of the twist
map $\gamma _{n}$ and put the signs of this map in a table. He partitioned
the twist table for \ $\mathbb{Z}_{2}^{n}$ into $2\times 2$ matrices and
obtained the following result:\medskip

\textbf{Theorem 1.2. (\cite{Ba-09}, Theorem 2.2., p. 88-91)} \textit{For} $%
n>0,$ \textit{the Cayley-Dickson twist table} $\gamma _{n}$ \textit{can be
partitioned in quadratic} \textit{matrices of dimension} $2$ \textit{of the
form} $A,B,C,-B,-C,$ \textit{defined in the relation} (\ref{1.6}). \textit{%
Relations between them can be found in the below twist trees}:

\begin{figure}[htbp]
\begin{center}
\includegraphics[height=7cm ,width=4cm, angle=270]{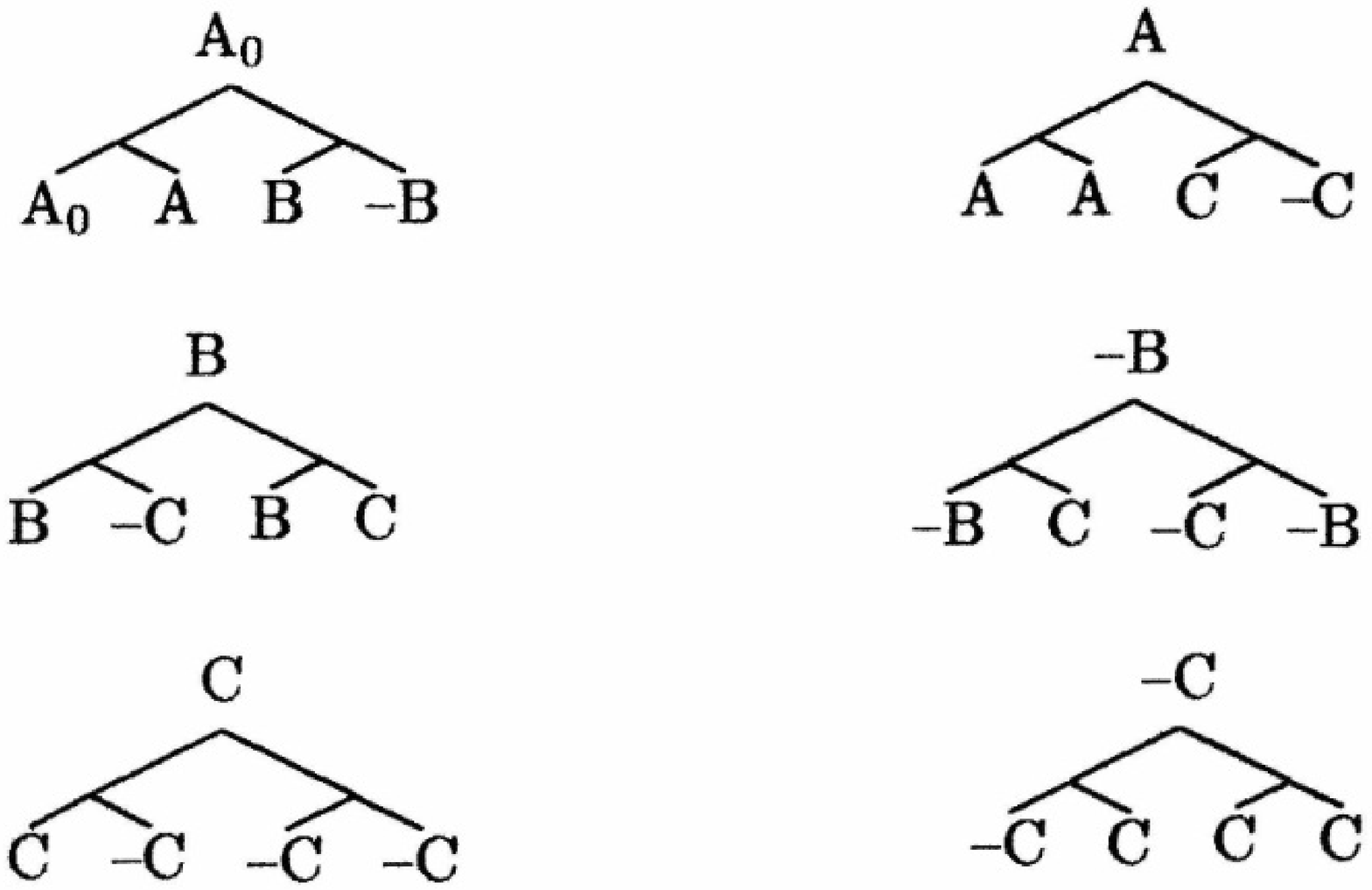}
\par
Fig. 1: Twist trees([1], Table 9)
\end{center}
\end{figure}

$\Box \medskip $

\textbf{Definition 1.3.} Let $x=x_{0},x_{1},x_{2},....$ and~ $%
y=y_{0},y_{1},y_{2},.....$ be two sequences of real numbers. The ordered
pair 
\begin{equation*}
\left( x,y\right) =x_{0},y_{0},x_{1},y_{1},x_{2},y_{2},....
\end{equation*}%
is a sequence obtained by \textit{shuffling} \ the sequences $x$ and$~\ y.$

Using Theorem 1.2, in \cite{Ba-09}, the author gave the below algorithm for
find $\gamma _{n}\left( s,r\right) ,$ where $s,r\in $ $\mathbb{Z}_{2}^{n}$ $%
: $

i) We find the shuffling sequence $\left( s,r\right).$

ii) Starting with the root $A_{0},$ we can find $\gamma _{n}\left(
s,r\right) $ using the twist tree. We remark that $"00"$= unchanged, $"01"$
=left $\rightarrow $right, $"10"$=right$\rightarrow $ left, $"11"$=right $%
\rightarrow $ right. \vspace{2mm}

Let $\mathbb{H}\left( \gamma _{1},\gamma _{2}\right) $ be the generalized
quaternion{\small \ }algebra and $\mathbb{H(}-1,-1)$ be the quaternion
division algebra. Below, you can see the multiplication tables:\medskip\ 
\begin{eqnarray*}
&&%
\begin{tabular}{l|llll}
$\cdot $ & $1$ & $e_{1}$ & $e_{2}$ & $e_{3}$ \\ \hline
$1$ & $1$ & $e_{1}$ & \multicolumn{1}{|l}{$e_{2}$} & \multicolumn{1}{l|}{$%
e_{3}$} \\ 
$e_{1}$ & $e_{1}$ & $\gamma _{1}$ & \multicolumn{1}{|l}{$e_{3}$} & 
\multicolumn{1}{l|}{$\gamma _{1}e_{2}$} \\ \cline{2-5}
$e_{2}$ & $e_{2}$ & $-e_{3}$ & \multicolumn{1}{|l}{$\gamma _{2}$} & 
\multicolumn{1}{l|}{$-\gamma _{2}e_{1}$} \\ 
$e_{3}$ & $e_{3}$ & $-\gamma _{1}e_{2}$ & \multicolumn{1}{|l}{$\gamma
_{2}e_{1}$} & \multicolumn{1}{l|}{$\gamma _{1}\gamma _{2}$} \\ \cline{2-5}
\end{tabular}%
\ \ \  \\
&&\text{{\small Multiplication} ~{\small table for the generalized
~quaternion~ algebra}}{\small \ }\ 
\end{eqnarray*}%
\begin{eqnarray*}
&&\ 
\begin{tabular}{l|llll}
$\cdot $ & $1$ & $e_{1}$ & $e_{2}$ & $e_{3}$ \\ \hline
$1$ & $1$ & $e_{1}$ & \multicolumn{1}{|l}{$e_{2}$} & \multicolumn{1}{l|}{$%
e_{3}$} \\ 
$e_{1}$ & $e_{1}$ & $-1$ & \multicolumn{1}{|l}{$e_{3}$} & 
\multicolumn{1}{l|}{$-e_{2}$} \\ \cline{2-5}
$e_{2}$ & $e_{2}$ & $-e_{3}$ & \multicolumn{1}{|l}{$-1$} & 
\multicolumn{1}{l|}{$e_{1}$} \\ 
$e_{3}$ & $e_{3}$ & $e_{2}$ & \multicolumn{1}{|l}{$-e_{1}$} & 
\multicolumn{1}{l|}{$-1$} \\ \cline{2-5}
\end{tabular}%
\ \ \  \\
&&\text{{\small Multiplication~ table} {\small for} {\small the} {\small %
real division ~quaternion~ algebra}}{\small \ }
\end{eqnarray*}%
\begin{eqnarray*}
&&\ \ 
\begin{tabular}{l|llll}
$\cdot $ & $1$ & $e_{1}$ & $e_{2}$ & $e_{3}$ \\ \hline
$1$ & $1$ & $1$ & \multicolumn{1}{|l}{$1$} & \multicolumn{1}{l|}{$1$} \\ 
$e_{1}$ & $1$ & $-1$ & \multicolumn{1}{|l}{$1$} & \multicolumn{1}{l|}{$-1$}
\\ \cline{2-5}
$e_{2}$ & $1$ & $-1$ & \multicolumn{1}{|l}{$-1$} & \multicolumn{1}{l|}{$1$}
\\ 
$e_{3}$ & $1$ & $1$ & \multicolumn{1}{|l}{$-1$} & \multicolumn{1}{l|}{$-1$}
\\ \cline{2-5}
\end{tabular}%
\ \ \ {\small \ } \\
&&\text{{\small Quaternion\ \ twist\ \ table}}
\end{eqnarray*}%
\begin{eqnarray*}
&&\left( 
\begin{array}{cc}
A_{0} & A \\ 
B & -B%
\end{array}%
\right) \ \  \\
&&\text{{\small Quaternion\ \ twist\ \ table using notations from Theorem
1.2.}}
\end{eqnarray*}

\bigskip

\textbf{Example 1.4.} Let $A_{4}$ be the real sedenion algebra. That means $%
\dim A_{4}=16$ with $\{1,e_{1},...,e_{15}\}$ a basis in this algebra. Let
compute $e_{7}e_{13}=\gamma _{4}(7_{2},13_{2})e_{7\otimes 13}.$ We have the
following binary decompositions: 
\begin{eqnarray*}
7_{2} &=&0111,\text{ since }7=2^{2}+2+1\text{ and } \\
13_{2} &=&1101,\text{ since }13=2^{3}+2^{2}+1.
\end{eqnarray*}

Since $0111\otimes 1101=1010 (=2^{3}+2=10),$ it results that $7\otimes
13=10. $

Now, we compute $\gamma _{4}\left( e_{7},e_{13}\right) .$ First, we shuffle
the sequences $0111$ and $1101.$ We obtain $01~11~10~11.$ Starting with $%
A_{0},$ it results: $A_{0}\overset{01}{\rightarrow }A\overset{11}{
\rightarrow }-C\overset{10}{\rightarrow }C\overset{11}{\rightarrow }-C,$
then $\gamma _{4}\left( e_{7},e_{13}\right) =-1$ and $e_{7}e_{13}=-e_{10}. $

\vspace{5mm}%
\begin{equation*}
\end{equation*}
\begin{equation*}
\end{equation*}
\begin{equation*}
\end{equation*}

\textbf{2. Main results}

\bigskip 
\begin{equation*}
\end{equation*}

In this section, for a generalized Cayley-Dickson algebra $A_{t},$\ writing
the basis's elements in a convenient way, we can obtain multiplication
tables for certain elements of \ the basis. \ Using these results, in
Theorem 2.12 \ we provide an example of a left\textit{\ }hyperholomorphic
function \ in generalized \ Cayley-Dickson algebras.$\medskip $

\bigskip \textbf{Remark 2.1.} i) In the generalized quaternion algebra, $%
\mathbb{H}\left( \gamma _{1},\gamma _{2}\right) ,$ the basis can be written
as 
\begin{equation*}
\{1=e_{0},e_{1},e_{2},e_{1}e_{2}\}.
\end{equation*}%
For the generalized octonion algebra, $\mathbb{O(}\gamma _{1},\gamma
_{2},\gamma _{3}),$ the basis can be written 
\begin{equation*}
\{1=e_{0},e_{1},e_{2},e_{1}e_{2},e_{4},e_{1}e_{4},e_{2}e_{4},\left(
e_{1}e_{2}\right) e_{4}\}.
\end{equation*}%
Therefore $e_{3}=e_{1}e_{2},e_{7}=e_{3}e_{4}=\left( e_{1}e_{2}\right) e_{4},$
$e_{2}e_{4}=e_{6}$ and, when compute them, in these products do not appear
any of the elements $\gamma _{1},\gamma _{2},\gamma _{3},$ or products of
some of them at the end.

We remark that in the algebra $A_{t}=\left( \frac{\gamma _{1},...,\gamma _{t}%
}{\mathbb{R}}\right) $ in the products of the form 
\begin{equation*}
e_{1}e_{2},\left( e_{1}e_{2}\right) e_{4},...,((e_{2^{r}}e_{2^{r+1}})\ldots
e_{2^{k}})e_{2^{i}},
\end{equation*}%
when compute them, do not appear any of the elements $\gamma _{1},\gamma
_{2},...,\gamma _{t}$ or products of some of them at the end.

ii) Using above remarks, the basis in the algebra $A_{t}=\left( \frac{\gamma
_{1},\ldots ,\gamma _{t}}{\mathbb{R}}\right) $ can be written under the form%
\begin{equation}
\{1=e_{0},e_{1},e_{2},\ldots
,e_{2^{t-1}-1},e_{2^{t-1}},e_{1}e_{2^{t-1}},e_{2}e_{2^{t-1}},e_{3}e_{2^{t-1}},\ldots ,e_{2^{t-1}-1}e_{2^{t-1}}\}
\label{2.1}
\end{equation}%
with 
\begin{equation}
e_{i}e_{2^{t-1}}=-e_{2^{t-1}}e_{i}=e_{2^{t-1}}\overline{e}_{i},\quad i\in
\{1,2,\ldots ,2^{t-1}-1\}.  \label{2.2}
\end{equation}

\textbf{Proposition 2.2. } \textit{Let} $A_{t}=\left( \frac{\gamma
_{1},...,\gamma _{t}}{\mathbb{R}}\right) $ \textit{be an algebra obtained by
the Cayley-Dickson process and} $\{e_{0}=1,e_{1},...,e_{n-1}\},\,\,n=2^{t}$ 
\textit{be a basis. Let} $\ r\geq 1,\,\,r<k\leq i<t.$ \textit{Therefore} 
\begin{equation}  \label{2.3}
((e_{2^{r}}e_{2^{r+1}})\ldots e_{2^{k}})e_{2^{i}}=\left( -1\right)
^{k-r+2}e_{T},
\end{equation}
\begin{equation}  \label{2.4}
((e_{1}e_{2^{r}})e_{2^{r+1}})\ldots e_{2^{k}})e_{2^{i}}=\left( -1\right)
^{k-r+3}e_{T+1},
\end{equation}
\textit{where} $T=2^{r}+2^{r+1}+\ldots +2^{k}+2^{i}$ \textit{and} 
\begin{equation}  \label{2.5}
e_{1}e_{2^{i}}=e_{2^{i}+1}.
\end{equation}
$.\medskip $

\textbf{Proof.} \ From Remark 2.1, it results that we can use Theorem 1.2 \
for $\gamma _{1},\gamma _{2},\ldots ,\gamma _{t}$ arbitrary. From Remark
1.1, it results $T=2^{r}+2^{r+1}+\ldots +2^{k}+2^{i}.$ For $T,$ we have the
binary decomposition 
\begin{equation*}
T_{2}=1\underset{i-k-1}{\underbrace{00\ldots 0}}\underset{k-r+1}{\underbrace{%
111\ldots 1}}\underset{r}{\underbrace{0\ldots 0}}.
\end{equation*}%
Using the same remark, we obtain $e_{2^{r}}e_{2^{r+1}}=\gamma _{n}\left( 
\underset{r+2}{\underbrace{01...0}},\underset{r+2}{\underbrace{10\ldots 0}}%
\right) e_{2^{r}+2^{r+1}}.$ We "shuffling" $\underset{r+2}{\underbrace{%
01\ldots 0}}\ \ $and\ $\underset{r+2}{\underbrace{10\ldots 0}}$ and we
obtain $01~10~\underset{r\,\,\text{pairs}}{\underbrace{00\ 00\ldots 00\ 00}}%
. $ Starting with $A_{0},$ it results: 
\begin{equation*}
A_{0}\overset{01}{\rightarrow }A\overset{10}{\rightarrow }C,
\end{equation*}%
then $\gamma _{n}\left( \underset{r+2}{\underbrace{01\ldots 0}},\underset{r+2%
}{\underbrace{10\ldots 0}}\right) =1$ and $%
e_{2^{r}}e_{2^{r+1}}=e_{2^{r}+2^{r+1}}.$

We compute $(e_{2^{r}}e_{2^{r+1}})e_{2^{r+2}}.$ We obtain 
\begin{equation*}
(e_{2^{r}}e_{2^{r+1}})e_{2^{r+2}}=e_{2^{r}+2^{r+1}}e_{2^{r+2}}=\gamma
_{n}\left( \underset{r+3}{\underbrace{011...0}},\underset{r+3}{\underbrace{%
10...0}}\right) e_{2^{r}+2^{r+1}+2^{r+2}}.
\end{equation*}%
Shuffling $\underset{r+3}{\underbrace{011...0}}$ and $\underset{r+3}{%
\underbrace{10...0}},$ we get $01~10~10\underset{r\,\,\text{pairs}}{%
\underbrace{00\ 00...00\ 00}}.$ Starting with $A_{0},$ it results: $A_{0}%
\overset{01}{\rightarrow }A\overset{10}{\rightarrow }C\overset{10}{%
\rightarrow }-C,$ then 
\begin{equation*}
\gamma _{n}\left( \underset{r+3}{\underbrace{011...0}},\underset{r+3}{%
\underbrace{10...0}}\right) =-1,
\end{equation*}%
therefore $e_{2^{r}+2^{r+1}}e_{2^{r+2}}=-e_{2^{r}+2^{r+1}+2^{r+2}}.$
Continuing this procedure, we remark that the number of "$1$" in the
"shuffling" obtained influences the sign. Since $%
T=2^{r}+2^{r+1}+...+2^{k}+2^{i}$ has binary decomposition 
\begin{equation*}
T_{2}=1\underset{i-k-1}{\underbrace{00...0}}\underset{k-r+1}{\underbrace{%
111..1}}\underset{r}{\underbrace{0...0}},
\end{equation*}%
in which we have $k-r+2$ elements equal with $1,$ we obtain relation (\ref%
{2.3}). In the same way it results relations (\ref{2.4}) and (\ref{2.5}). $%
\Box \medskip $

\textbf{Proposition 2.3.} \textit{With the same notations as in Proposition
2.2, for the algebra }$A_{t}=\left( \frac{-1,...,-1}{\mathbb{R}}\right) ,$%
\thinspace\ \textit{we have:} 
\begin{equation}
\begin{tabular}{l|ll}
$\cdot $ & $e_{T}$ & $e_{T+1}$ \\ \hline
$e_{T_{1}}$ & $\left( -1\right) ^{k-r+1}e_{2^{i}}$ & $-\left( -1\right)
^{k-r+1}e_{2^{i}+1}$ \\ 
$e_{T_{1}+1}$ & $-\left( -1\right) ^{^{k-r+1}}e_{2^{i}+1}$ & $-\left(
-1\right) ^{^{k-r+1}}e_{2^{i}}$%
\end{tabular}
\label{2.6}
\end{equation}%
\textit{for} $r<k$, \textit{where} $T=2^{r}+2^{r+1}+...+2^{k}+2^{i},$ $%
T_{1}=2^{r}+2^{r+1}+...+2^{k}$ \textit{and} 
\begin{equation}
\begin{tabular}{l|ll}
$\cdot $ & $e_{T}$ & $e_{T+1}$ \\ \hline
$e_{2^{k}}$ & $e_{M}$ & $-e_{M+1}$ \\ 
$e_{2^{k}+1}$ & $-e_{M+1}$ & $-e_{M}$%
\end{tabular}%
,  \label{2.7}
\end{equation}%
\textit{where} $M=2^{k}+2^{i}\medskip .\medskip $

\textbf{Proof. }\textit{Case 1}: $r<k.$ We \ compute $e_{T_{1}}e_{T}.$ We
have $e_{T_{1}}e_{T}=\gamma \left( s,q\right) e_{M},$ where $s,q$ are the
binary decomposition of $T_{1}$ and $T.$ The binary decomposition of $M$ is$%
\ M_{2}=T_{1}\otimes T.$ It results $M$ $=2^{i},$ 
\begin{equation*}
s=~\underset{i-k}{\underbrace{00...0}}\underset{k-r+1}{\underbrace{111...1}}%
\underset{r}{\underbrace{0...0}},\quad q=\ \underset{i-k}{\underbrace{100...0%
}}\underset{k-r+1}{\underbrace{111...1}}\underset{r}{\underbrace{0...0}}.
\end{equation*}%
By "shuffling" $s\otimes q$, we obtain 
\begin{equation*}
\underset{(i-k)\,\,\text{pairs}}{\underbrace{01~00~00...00}\ }\underset{%
(k-r+1)\,\,\text{pairs}}{\ \underbrace{11~11~11~...11}}\underset{r\,\,\text{%
pairs}}{\ \underbrace{00~00~...00~00}}.
\end{equation*}

Starting with $A_{0},$ we get:

\begin{equation*}
\underset{i-k}{\underbrace{A_{0}\overset{01}{\rightarrow }A\overset{00}{
\rightarrow }...\overset{00}{\rightarrow }}}\underset{k-r+1}{\underbrace{A 
\overset{11}{\rightarrow }-C\overset{11}{\rightarrow }C\overset{11}{
\rightarrow }-C\overset{11}{\rightarrow }C\overset{11}{\rightarrow }... 
\overset{11}{\rightarrow } \left( -1\right) ^{k-r+1}C}}\underset{r}{%
\underbrace{\overset{00}{ \rightarrow }...\overset{00}{\rightarrow }\left(
-1\right) ^{k-r+1}C}}.
\end{equation*}
Therefore $\gamma \left( s,q\right) =\left( -1\right) ^{k-r+1}.$

Now, we compute $e_{T_{1}}e_{T+1}.$ For this, we will "shuffling" $\underset{
i-k}{\underbrace{00...0}}\underset{k-r+1}{\underbrace{111...1}}\underset{r}{ 
\underbrace{0...0}}$ with $\underset{i-k}{\underbrace{100...0}}\underset{
k-r+1}{\underbrace{111...1}}\underset{r}{\underbrace{0...1}}.$ It results 
\begin{equation*}
\underset{(i-k)\,\,\text{pairs}}{\underbrace{01~00~00...00}\ }\underset{%
(k-r+1)\,\,\text{pairs} }{\ \underbrace{ 11~11~11...11}}\underset{r \,\,%
\text{pairs}}{\ \underbrace{00~00...00~01}}.
\end{equation*}
Starting with $A_{0},$ we get:

\begin{equation*}
\underset{i-k}{\underbrace{A_{0}\overset{01}{\rightarrow }A\overset{00}{
\rightarrow }...\overset{00}{\rightarrow }}}\underset{k-r+1}{\underbrace{A 
\overset{11}{\rightarrow }-C\overset{11}{\rightarrow }C\overset{11}{
\rightarrow }-C\overset{11}{\rightarrow }C\overset{11}{\rightarrow } ...%
\overset{11}{\rightarrow } \left( -1\right) ^{k-r+1}C}}\underset{r}{%
\underbrace{\overset{00}{ \rightarrow }...\overset{01}{\rightarrow }-\left(
-1\right) ^{k-r+1}C}}.
\end{equation*}

For $e_{T_{1}+1}e_{T},$ "shuffling" $\underset{i-k}{\underbrace{00...0}} 
\underset{k-r+1}{\underbrace{111...1}}\underset{r}{\underbrace{0...1}}$ with 
$\underset{i-k}{\underbrace{100...0}}\underset{k-r+1}{\underbrace{111...1}} 
\underset{r}{\underbrace{0...0}},$ it results 
\begin{equation*}
\underset{(i-k)\,\,\text{pairs}}{\underbrace{01~00~00...00}\ }\underset{%
(k-r+1)\,\,\text{pairs} }{\ \underbrace{ 11~01~01...01}}\underset{r\,\,\text{%
pairs}}{\ \underbrace{00~00...00~10}}.
\end{equation*}
$\ $Starting with $A_{0},$we get: 
\begin{equation*}
\underset{i-k}{\underbrace{A_{0}\overset{01}{\rightarrow }A\overset{00}{
\rightarrow }...\overset{00}{\rightarrow }}}\underset{k-r+1}{\underbrace{A 
\overset{11}{\rightarrow }-C\overset{11}{\rightarrow }C\overset{11}{
\rightarrow }-C\overset{11}{\rightarrow }C\rightarrow...\overset{11}{%
\rightarrow } \left( -1\right) ^{k-r+1}C}}\underset{r}{\underbrace{\overset{%
00}{ \rightarrow }...\overset{10}{\rightarrow }-\left( -1\right) ^{k-r+1}C}}.
\end{equation*}
For $e_{T_{1}+1}e_{T+1},$ we compute first $\left( T_{1}+1\right) \otimes
\left( T+1\right).$ We obtain:\newline
\begin{equation*}
\left( 2^{r}+2^{r+1}+...+2^{k}+1\right) \otimes \left(
2^{r}+2^{r+1}+...+2^{k}+2^{i}+1\right) =\newline
\end{equation*}
\begin{equation*}
=\left( \underset{i-k}{\underbrace{00...0}}\underset{k-r+1}{\underbrace{
111...1}}\underset{r}{\underbrace{0...1}}\right) \otimes \left( \underset{i-k%
}{\underbrace{100...0}}\underset{k-r+1}{\underbrace{111...1}}\underset{r}{ 
\underbrace{0...1}}\right) =\medskip \medskip \newline
\end{equation*}
\begin{equation*}
=\ \ \ \underset{i-k}{\underbrace{10...0}}\underset{k-r+1}{\underbrace{
000...0}}\underset{r}{\underbrace{0...0}}\ =2^{i}.
\end{equation*}
$\medskip $

Now, "shuffling" $\underset{i-k}{\underbrace{00...0}}\underset{k-r+1}{ 
\underbrace{111...1}}\underset{r}{\underbrace{0...1}}$ with $\underset{i-k}{ 
\underbrace{100...0}}\underset{k-r+1}{\underbrace{111...1}}\underset{r}{ 
\underbrace{0...1}},$ it results 
\begin{equation*}
\underset{(i-k)\,\,\text{pairs} }{\underbrace{01~00~00...00}\ }\underset{%
(k-r+1)\,\,\text{pairs}}{\ \underbrace{ 11~01~01...01}}\underset{r\,\,\text{%
pairs}}{\ \underbrace{00~00...00~11}}.
\end{equation*}
\newline
Starting with $A_{0},$ we get:

\begin{equation*}
\underset{i-k}{\underbrace{A_{0}\overset{01}{\rightarrow }A\overset{00}{%
\rightarrow }...\overset{00}{\rightarrow }}}\underset{k-r+1}{\underbrace{A%
\overset{11}{\rightarrow }-C\overset{11}{\rightarrow }C\overset{11}{%
\rightarrow }-C\overset{11}{\rightarrow }C\overset{11}{\rightarrow }...%
\overset{11}{\rightarrow }\left( -1\right) ^{k-r+1}C}}\underset{r}{%
\underbrace{\overset{00}{\rightarrow }...\overset{11}{\rightarrow }-\left(
-1\right) ^{k-r+1}C}}.
\end{equation*}%
\textit{Case 2}: $r=k.~$We have $M=2^{k}\otimes T=2^{i}+2^{k}.$ For $%
e_{2^{k}}e_{T},$ "shuffling" $\underset{i-k}{\underbrace{00...0}}\underset{%
k+1}{\underbrace{10...0}}$ with $\underset{i-k}{\underbrace{100...0}}%
\underset{k+1}{\underbrace{0...0}},$ it results 
\begin{equation*}
\underset{(i-k)\,\,\text{pairs}}{\underbrace{01~00~00...00}\ }\underset{%
(k+1)\,\,\text{pairs}}{\ \underbrace{10~00~00~...00}}.
\end{equation*}%
$\ $ Starting with $A_{0},$ we get: 
\begin{equation*}
\underset{i-k}{\underbrace{A_{0}\overset{01}{\rightarrow }A\overset{00}{%
\rightarrow }.....\overset{00}{\rightarrow }}}\underset{k+1}{\underbrace{A%
\overset{10}{\rightarrow }C\overset{00}{\rightarrow }C\overset{00}{%
\rightarrow }...\overset{00}{\rightarrow }C}}.
\end{equation*}%
For $e_{2^{k}}e_{T+1},$ "shuffling" $\underset{i-k}{\underbrace{00...0}}%
\underset{k+1}{\underbrace{10...0}}$ with $\underset{i-k}{\underbrace{100...0%
}}\underset{k+1}{\underbrace{0...1}},$ it results 
\begin{equation*}
\underset{(i-k)\,\,\text{pairs}}{\underbrace{01~00~00...00}\ }\underset{%
(k+1)\,\,\text{pairs}}{\ \underbrace{10~00~00~...01}}.
\end{equation*}%
$\ $Starting with $A_{0},$ we get: 
\begin{equation*}
\underset{i-k}{\underbrace{A_{0}\overset{01}{\rightarrow }A\overset{00}{%
\rightarrow }}...\overset{00}{\rightarrow }}\underset{k+1}{\underbrace{A%
\overset{10}{\rightarrow }C\overset{00}{\rightarrow }C\overset{00}{%
\rightarrow }...\overset{01}{\rightarrow }-C}}.
\end{equation*}%
etc.$\Box \medskip $

\textbf{Proposition 2.4.} \textit{Let} $A_{t}=\left( \frac{\gamma
_{1},...,\gamma _{t}}{\mathbb{R}}\right) $ \textit{be an algebra obtained by
the Cayley-Dickson process. For any }$x_{1},x_{2},...,x_{t}\in \mathbb{R}%
-\{0\},$ \textit{we have that} 
\begin{equation*}
\left( \frac{\gamma _{1},...,\gamma _{t}}{\mathbb{R}}\right) \simeq \left( 
\frac{\gamma _{1}x_{1}^{2},...,\gamma _{t}x_{t}^{2}}{\mathbb{R}}\right) .
\end{equation*}

\textbf{Proof.} Let $A_{t}=\left( \frac{\gamma _{1},...,\gamma _{t}}{\mathbb{%
\ R}}\right) $ with the basis $\{e_{0}=1,e_{1},...,e_{n-1}\},\,\,n=2^{t}$
and let $A_{t}^{\prime }=\left( \frac{\gamma _{1}x_{1}^{2},...,\gamma
_{t}x_{t}^{2}}{\mathbb{R}}\right) $ with the basis $\{e_{0}^{\prime
}=1,e_{1}^{\prime },...,e_{n-1}^{\prime }\}$ such that $(e_{i}^{\prime
})^{2}=\gamma _{i}x_{i}^{2},\,\,i\in \{1,2,...,n-1\}.$ We remark that $%
(x_{i}e_{i})^{2}=x_{i}^{2}\gamma _{i}$ and from here, it results that the
map $\tau :A_{t}^{\prime }\rightarrow A_{t},\,\,\tau \left( e_{i}^{\prime
}\right) =e_{i}x_{i}$ is \ an algebra isomorphism. $\Box \medskip $

The above proposition generalized Proposition 1.1, p. 52 from \cite{La-05}%
.\medskip

\textbf{Remark 2.5.} From Proposition 2.4, it results that for each $n=2^{t}$
there are only $n$ non-isomorphic algebras $A_{t}$. These algebras are of
the form $A_{t}=\left( \frac{\gamma _{1},...,\gamma _{t}}{\mathbb{R}}\right) 
$, with $\gamma _{1},...,\gamma _{t}\in \{-1,1\}.\medskip $

\textbf{Definition 2.6.} Let $\{e_{0}=1,e_{1},...,e_{n-1}\}$ be a basis in $%
A_{t}=\left( \frac{\gamma _{1},...,\gamma _{t}}{\mathbb{R}}\right) ,$
\thinspace\ $n=2^{t}$. To domain $\Omega \subset \mathbb{R}^{2^{t}-1}$ we
will associate the domain $\Omega _{\zeta }:=\{\zeta =x_{1}e_{1}+\ldots
+x_{n-1}e_{n-1}:(x_{1},x_{2},\ldots ,x_{n-1})\in \Omega \}$ included in $%
A_{t}.$

Consider a function $\Phi:\Omega_{\zeta}\rightarrow A_t$ of the form 
\begin{equation}  \label{Fi}
\Phi(\zeta)=\sum\limits_{k=1}^{n-1}\Phi_{k}(x_1,x_2,\ldots,x_{n-1})e_{k},
\end{equation}
where $(x_1,x_2,\ldots, x_{n-1})\in \Omega$ and $\Phi_k : \Omega\rightarrow%
\mathbb{R}$.

We say that a function of the form (\ref{Fi}) is \textit{left} $A_{t}-$%
\textit{holomorphic} in a domain $\Omega _{\zeta }$ if the first partial
derivatives $\partial \Phi _{k}/\partial x_{k}$ exist in $\Omega $ and the
following equality is fulfilled in every point of $\Omega _{\zeta }$: 
\begin{equation*}
D[\Phi ](\zeta )=\underset{k=1}{\overset{2^{t}-1}{\sum }}e_{k}\frac{\partial
\Phi }{\partial x_{k}}=0.
\end{equation*}%
The operator $D$ is called \textit{Dirac operator}. Note that if $A_{t}$ is
the generalized quaternion algebra, then the left $A_{t}-$holomorphic
functions is also called \textit{hyperholomorphic}. We also note that every
hyperholomorphic function $\Phi $ in a domain $\Omega _{\zeta }$ is a
solution of the equation 
\begin{equation*}
\gamma _{1}\frac{\partial ^{2}\Phi }{\partial x_{1}^{2}}+\gamma _{2}\frac{%
\partial ^{2}\Phi }{\partial x_{2}^{2}}+\gamma _{1}\gamma _{2}\frac{\partial
^{2}\Phi }{\partial x_{3}^{2}}=0.
\end{equation*}

\textbf{Remark 2.7.} Let $\mathbb{H}\left( \gamma _{1},\gamma _{2}\right) $
be the generalized quaternion{\small \ }algebra with the basis $%
\{1,e_{1},e_{2},e_{3}\},\,\,\gamma _{1}>0,\,\,\gamma _{2}>0\ $and $\mathbb{H(%
}-1,-1)$ be the usual quaternion division algebra with the basis $%
\{1,i,j,k\}.$ Let $\Omega $ be a domain in $\mathbb{R}^{3}$, and let $\Omega
_{\zeta }:=\{\zeta =xi+yj+zk:(x,y,z)\in \Omega \}$ be a corresponding domain
in $\mathbb{H(}-1,-1)$. The function $\Phi :\Omega _{\zeta }\rightarrow 
\mathbb{H(}-1,-1)$ of the form 
\begin{equation*}
\Phi (\zeta )=u_{1}\left( x,y,z\right) +u_{2}\left( x,y,z\right)
i+u_{3}\left( x,y,z\right) j+u_{4}\left( x,y,z\right) k.
\end{equation*}%
is hyperholomorphic in the domain $\Omega $ if 
\begin{equation*}
D[\Phi ](\zeta )=i\frac{\partial \Phi }{\partial x}+j\frac{\partial \Phi }{%
\partial y}+k\frac{\partial \Phi }{\partial z}=0.
\end{equation*}

For another domain $\Delta \subset \mathbb{R}^{3},$ we associate the domain $%
\Delta _{\widetilde{\zeta }}:=\{\widetilde{\zeta }=\widetilde{x}e_{1}+%
\widetilde{y}e_{2}+\widetilde{z}e_{3}:(\widetilde{x},\widetilde{y},%
\widetilde{z})\in \Delta \}$ in the algebra $\mathbb{H}\left( \gamma
_{1},\gamma _{2}\right) $. The Dirac operator in $\mathbb{H}\left( \gamma
_{1},\gamma _{2}\right) ,$ denoted by $\widetilde{D}$, is 
\begin{equation*}
\widetilde{D}:=e_{1}\frac{\partial }{\partial \widetilde{x}}+e_{2}\frac{%
\partial }{\partial \widetilde{y}}+e_{3}\frac{\partial }{\partial \widetilde{%
z}}.
\end{equation*}

The elements of bases in $\mathbb{H}\left( -1,-1\right) $ and $\mathbb{H}%
\left( \gamma _{1},\gamma _{2}\right) $ satisfy the following equalities: 
\begin{equation}
e_{1}=i\sqrt{\gamma _{1}},\quad e_{2}=j\sqrt{\gamma _{2}},\quad e_{3}=k\sqrt{%
\gamma _{1}\gamma _{2}}.  \label{2.8}
\end{equation}

Now we establish a connection between hyperholomorphic functions in the
algebras $\mathbb{H}\left( -1,-1\right) $ and $\mathbb{H}\left( \gamma
_{1},\gamma _{2}\right) $, where $\gamma _{1}>0,\,\,\gamma _{2}>0$. For
this, we denote 
\begin{equation*}
x=\frac{1}{\sqrt{\gamma _{1}}}\widetilde{x},\quad y=\frac{1}{\sqrt{\gamma
_{2}}}\widetilde{y},\quad z=\frac{1}{\sqrt{\gamma _{1}\gamma _{2}}}%
\widetilde{z}.
\end{equation*}

These relations give us the operator equalities: 
\begin{equation}
\frac{\partial }{\partial \widetilde{x}}=\frac{1}{\sqrt{\gamma _{1}}}\frac{%
\partial }{\partial x},\quad \frac{\partial }{\partial \widetilde{y}}=\frac{1%
}{\sqrt{\gamma _{2}}}\frac{\partial }{\partial y},\quad \frac{\partial }{%
\partial \widetilde{z}}=\frac{1}{\sqrt{\gamma _{1}\gamma _{2}}}\frac{%
\partial }{\partial z}.  \label{2.9}
\end{equation}

Now, using relations (\ref{2.8}) and (\ref{2.9}), we obtain 
\begin{equation*}
\widetilde{D}[\Phi](\widetilde{\zeta}\,) =e_{1}\frac{\partial \Phi }{
\partial \widetilde{x}}+e_{2}\frac{\partial \Phi }{\partial \widetilde{y}}
+e_{3}\frac{\partial \Phi }{\partial \widetilde{z}}=
\end{equation*}
\begin{equation*}
=i\frac{\partial \Phi }{\partial x}\frac{1}{\sqrt{\gamma _{1}}}\sqrt{\gamma
_{1}}+j\frac{\partial \Phi }{\partial y}\frac{1}{\sqrt{\gamma _{2}}}\sqrt{
\gamma _{2}}+k\frac{\partial \Phi }{\partial z}\frac{1}{\sqrt{\gamma
_{1}\gamma _{2}}}\sqrt{\gamma _{1}\gamma _{2}}=
\end{equation*}
\begin{equation*}
=i\frac{\partial \Phi }{\partial x}+j\frac{\partial \Phi }{\partial y}+k%
\frac{ \partial \Phi }{\partial z}=D[\Phi](\zeta)=0.
\end{equation*}

Using the above notations, we obtain the following theorem: \medskip

\textbf{Theorem 2.8.} \textit{Let }$\Omega $\textit{\ be an arbitrary domain
in }$\mathbb{R}^{3}$\textit{\ and }$\Delta $\textit{\ be a domain in }$%
\mathbb{R}^{3}$\textit{\ such that the coordinates of the corresponding
points }$\zeta =xi+yj+zk\in \Omega _{\zeta }$\textit{\ and }$\widetilde{%
\zeta }=\widetilde{x}e_{1}+\widetilde{y}e_{2}+\widetilde{z}e_{3}\in \Delta _{%
\widetilde{\zeta }}$\textit{\ satisfy the following relations}:%
\begin{equation*}
x=\frac{1}{\sqrt{sign(\gamma _{1})\gamma _{1}}}\widetilde{x},\,\,\,y=\frac{1%
}{\sqrt{sign(\gamma _{2})\gamma _{2}}}\widetilde{y},\,\,\,z=\frac{1}{\sqrt{%
sign(\gamma _{1})sign(\gamma _{2})\gamma _{1}\gamma _{2}}}\widetilde{z},
\end{equation*}%
\textit{\ where }$sign(a)$\textit{\ is the sign of the non-zero real number }%
$a$\textit{$.$ Then if the function} $\Phi :\Omega _{\zeta }\rightarrow 
\mathbb{H}\big(sign(\gamma _{1}),sign(\gamma _{2})\big)$ \textit{is
hyperholomorphic in the domain }$\Omega _{\zeta }$\textit{, then the same
function }$\Phi $\textit{, of }$\widetilde{\zeta },$ \textit{is
hyperholomorphic in the domain }$\Delta _{\widetilde{\zeta }}\in \mathbb{H}%
(\gamma _{1},\gamma _{2})$\textit{. The converse is also true.\medskip }

\textbf{Proof.} Since $e_{1}=i\sqrt{sign(\gamma _{1})\gamma _{1}},$
\thinspace \thinspace $e_{2}=j\sqrt{sign(\gamma _{2})\gamma _{2}}$%
,\thinspace \thinspace\ \newline
$e_{3}=k\sqrt{sign(\gamma _{1})sign(\gamma _{2})\gamma _{1}\gamma _{2}},$
the result directly follows from Remark 2.7.$\Box \medskip \medskip $

\textbf{Remark 2.9. }(i)\thinspace The above Theorem tell us that for
studying hyperholomorphic functions in generalized quaternion algebras $%
\mathbb{H(}\gamma _{1},\gamma _{2})$ it is suffices to consider
hyperholomorphic functions only in the algebras $\mathbb{H}\big(sign(\gamma
_{1}),sign(\gamma _{2})\big).$

(ii)\thinspace \thinspace\ The result similar to the previous remark was
established in the paper \cite{Pl-Shp-11} (Theorem 5) in a three-dimensional
commutative associative algebra. \medskip

\textbf{Theorem 2.10.} \textit{Let} $A_{t}=\left( \frac{\gamma
_{1},...,\gamma _{t}}{\mathbb{R}}\right) $ \textit{be a generalized
Cayley-Dickson algebra.} \textit{Let }$\Omega $\textit{\ be an arbitrary
domain in }$\mathbb{R}^{2^{t}-1}$\textit{\ and }$\Delta $\textit{\ be a
domain in }$\mathbb{R}^{2^{t}-1}$\textit{\ such that the coordinates of the
corresponding points }$\zeta =x_{1}e_{1}+\ldots +x_{2^{t}-1}e_{2^{t}-1}\in
\Omega _{\zeta }$\textit{\ and }$\widetilde{\zeta }=\widetilde{x}_{1}\,%
\widetilde{e}_{1}+\widetilde{x}_{2}\,\widetilde{e}_{2}+\ldots +\widetilde{x}%
_{2^{t}-1}\,\widetilde{e}_{2^{t}-1}\in \Delta _{\widetilde{\zeta }}$\textit{%
\ satisfy the following relations}%
\begin{equation*}
x_{1}=\frac{1}{\sqrt{sign(\gamma _{1})\gamma _{1}}}\widetilde{x}_{1},\quad
x_{2}=\frac{1}{\sqrt{sign(\gamma _{2})\gamma _{2}}}\widetilde{x}_{2},\ldots
\end{equation*}%
\begin{equation*}
\ldots ,x_{n}=\frac{1}{\sqrt{sign(\gamma _{1})...sign(\gamma _{t})\gamma
_{1}...\gamma _{t}}}\widetilde{x}_{n}.
\end{equation*}%
\textit{\ If the function }$\Phi :\Omega _{\zeta }\rightarrow \left( \frac{%
sign(\gamma _{1}),...,sign(\gamma _{t})}{\mathbb{R}}\right) $\textit{\ is
left }$A_{t}$\textit{-holomorphic in the domain }$\Omega _{\zeta }$\textit{,
then the same function }$\Phi $\textit{, but depending of }$\widetilde{\zeta 
}$\textit{\ is left }$A_{t}$\textit{-holomorphic in the domain }$\Delta _{%
\widetilde{\zeta }}\in A_{t}$\textit{. The converse is also true.\medskip }

\textbf{Proof.} Let $\{1,e_{1},...,e_{n-1}\}$ be a basis in $\left( \frac{%
sign\left( \gamma _{1}\right),...,sign\left( \gamma _{t}\right) }{\mathbb{R} 
}\right)$ and $\{1,\widetilde{e}_1,...,\widetilde{e}_{n-1}\}$ be a basis in $%
A_{t}=\left( \frac{\gamma _{1},...,\gamma _{t}}{\mathbb{R}}\right) $.

Since 
\begin{equation*}
\widetilde{e}_{1}=e_{1}\sqrt{sign\left( \gamma _{1}\right) \gamma _{1}},
\quad \widetilde{e}_{2}=e_{2}\sqrt{sign\left( \gamma _{2}\right) \gamma _{2}}%
,\ldots,
\end{equation*}
\begin{equation*}
\ldots,\widetilde{e}_{n-1}=e_{n-1}\sqrt{sign\left( \gamma _{1}\right)
...sign\left( \gamma _{t}\right) \gamma _{1}...\gamma _{t}}\,,
\end{equation*}
the result is obtained from a simple computation as in Remark 2.7.$\Box
\medskip $

\textbf{Remark 2.11.} Using above Theorem, it is obvious that, for studying
left $A_{t}$-holomorphic functions in generalized Cayley-Dickson algebras $%
A_{t}=\left( \frac{\gamma _{1},...,\gamma _{t}}{\mathbb{R}}\right) $ it is
suffices to consider left $A_{t}$-holomorphic functions only in the algebras 
$\left( \frac{sign(\gamma _{1}),...,sign(\gamma _{t})}{\mathbb{R}}\right)
.\medskip $

Now we consider another class of differentiable functions. Let $A_{t}=\left( 
\frac{\gamma _{1},...,\gamma _{t}}{\mathbb{R}}\right) ,$ with $\gamma
_{1}=...=\gamma _{t}=-1,$ and the domain $\Omega \subset \mathbb{R}^{2^{t}}.$
We denote with $\Omega _{\zeta }:=\{\zeta =x_{0}+x_{1}e_{1}+\ldots
+x_{n-1}e_{n-1}:(x_{0},x_{1},\ldots ,x_{n-1})\in \Omega \}$ a domain in $%
A_{t}.$ This domain is called \textit{congruent} with the domain $\Omega .$

We consider a function $\Phi :\Omega _{\zeta }\rightarrow A_{t}$ of the form 
\begin{equation}
\Phi (\zeta )=\sum\limits_{k=0}^{n-1}\Phi _{k}(x_{0},x_{1},\ldots
,x_{n-1})e_{k},  \label{Fi-1}
\end{equation}%
where $(x_{0},x_{1},\ldots ,x_{n-1})\in \Omega $ and $\Phi _{k}:\Omega
\rightarrow \mathbb{R}$.

We say that a function of the form (\ref{Fi-1}) is \textit{left} $A_{t}-$%
\textit{hyperholomorphic} in a domain $\Omega _{\zeta }$ if the first
partial derivatives $\partial \Phi _{k}/\partial x_{k}$ exist in $\Omega $
and the following equality is fulfilled in every point of $\Omega _{\zeta }$%
: 
\begin{equation*}
\underset{k=0}{\overset{2^{t}-1}{\sum }}e_{k}\frac{\partial \Phi }{\partial
x_{k}}=0.
\end{equation*}

In the following, we will provide an algorithm to constructing a left $%
A_{t}- $hyperholomorphic functions. Using the above notations, let $v\left(
x,y\right) $ be a rational function defined in a domain $G\subset \mathbb{R}%
^{2}.$ In the following, using some ideas given in Theorem 3 from \cite%
{Xi-Zh-Li-05}, we will give an example of left $A_{t}-$hyperholomorphic
function, for all $t\geq 1,\,\,t\in \mathbb{N}$. For this, we consider the
following functions:%
\begin{equation*}
\phi _{1}=x_{0}+e_{1}x_{1},\quad \phi _{2}=\frac{1}{e_{1}}(x_{0}+e_{1}x_{1}),
\end{equation*}%
\begin{equation*}
\rho _{2s-1}=x_{2s}-e_{1}x_{2s+1},\quad \rho _{2s}=-\frac{1}{e_{1}}%
(x_{2s}-e_{1}x_{2s+1}),\,\,\,s\in \{1,2,...,2^{t-1}-1\},
\end{equation*}
\begin{equation*}
F_{t}\left( \zeta \right) =v\left( \phi _{1},\phi _{2}\right) +v\left( \rho
_{1},\rho _{2}\right) e_{2}+v\left( \rho _{3},\rho _{4}\right) e_{4}+\left[
v\left( \rho _{5},\rho _{6}\right) e_{2}\right] e_{4}+
\end{equation*}%
\begin{equation*}
+v\left( \rho _{7},\rho _{8}\right) e_{8}+\left( v\left( \rho _{9},\rho
_{10}\right) e_{2}\right) e_{8}+\left( v\left( \rho _{11},\rho _{12}\right)
e_{4}\right) e_{8}+\left[ \left( v\left( \rho _{13},\rho _{14}\right)
e_{2}\right) e_{4}\right] e_{8}+...
\end{equation*}%
\begin{equation*}
...+\underset{i=4}{\overset{t-1}{\sum }}\underset{k=1}{(\overset{i}{\sum }}%
\underset{r=1}{\overset{k-1}{(\sum }}v\left( \rho _{M_{rki}-1},\rho
_{M_{rki}}\right) e_{2^{r}})e_{2^{r+1}}...)e_{2^{k}})e_{2^{i}})+\underset{i=1%
}{\overset{t-1}{\sum }}\left( v\left( \rho _{2^{i}-1},\rho _{2^{i}}\right)
e_{2^{i}}\right) ,
\end{equation*}%
where $M_{rki}=2^{r}+2^{r+1}+...+2^{k}+2^{i}.\medskip $

It results

\begin{equation*}
F_{t}\left( \zeta \right) =v\left( \phi _{1},\phi _{2}\right) +
\end{equation*}%
\begin{equation*}
+\underset{i=1}{\overset{t-1}{\sum }}\underset{k=1}{(\overset{i}{\sum }}%
\underset{r=1}{\overset{k-1}{(\sum }}v\left( \rho _{M_{rki}-1},\rho
_{M_{rki}}\right) e_{2^{r}})e_{2^{r+1}}...)e_{2^{k}})e_{2^{i}})+\underset{i=1%
}{\overset{t-1}{\sum }}\left( v\left( \rho _{2^{i}-1},\rho _{2^{i}}\right)
e_{2^{i}}\right) ,
\end{equation*}

or

\begin{equation*}
F_{t}\left( \zeta \right) =F_{t-1}\left( \zeta \right) +
\end{equation*}%
\begin{equation*}
+\underset{k=1}{(\overset{t-2}{\sum }}\underset{r=1}{\overset{k-1}{(\sum }}%
v\left( \rho _{M_{rk(t-1)}-1},\rho _{M_{rk(t-1)}}\right)
e_{2^{r}})e_{2^{r+1}}...)e_{2^{k}})e_{2^{t-1}})+v\left( \rho
_{2^{t-1}-1},\rho _{2^{t-1}}\right) e_{2^{t-1}}.
\end{equation*}

\medskip

We denote with $\mathbb{C}_{2s}$ the "complex" planes $%
\{x_{2s}+e_{1}x_{2s+1}:x_{2s},x_{2s+1}\in \mathbb{R}\}$ and with $%
D_{2s}:=\{(x_{2s},x_{2s+1}):x_{2s}+e_{1}x_{2s+1}\in \mathbb{C}_{2s}\}$, $%
s\in \{0,1,2,...,2^{t-1}-1\}$ the Euclidian planes. Let $G_{2s}$ be a
domains in $\mathbb{C}_{2s}$ and let $\widetilde{G}_{2s}$ be the
corresponded domains in $D_{2s}$. We have the following theorem: \medskip

\textbf{Theorem 2.12.} \textit{With the above notations, we consider the
functions} $v\left( \phi _{1},\phi _{2}\right) $ \textit{and }$v\left( \rho
_{2s-1},\rho _{2s}\right) $ \textit{\ defined in the corresponding domains }$%
G_{0}\subset \mathbb{C}_{0}$\textit{\ and }$G_{2s}\subset \mathbb{C}_{2s},$ $%
s\in \{1,2,...,2^{t-1}-1\}$\textit{. Then the map }$F_{t}\left( \zeta
\right) $\textit{\ is a left }$A_{t}-$\textit{hyperholomorphic function in
the domain }$\Theta \subset A_{t}$\textit{\ which is congruent with the
domain }$\widetilde{G}_{0}\times \widetilde{G}_{2}\times \widetilde{G}%
_{4}\times ...\times \widetilde{G}_{2^{t-1}-1}\subset \mathbb{R}^{2^{t}},$%
\textit{\ for }$t\geq 1$\textit{$.$}

\medskip

\textbf{Proof.} For $t=1,$ we have $F_{1}\left( \zeta \right) =v\left( \phi
_{1},\phi _{2}\right) ,$ which is an holomorphic function in $D_{0}\subset 
\mathbb{C}_{0},$ as we can see in \cite{Xi-Zh-Li-05}, Theorem 3.

For $t=2,$ we obtain $F_{2}\left( \zeta \right) =v\left( \phi _{1},\phi
_{2}\right) +v\left( \rho _{1},\rho _{2}\right) e_{2}$ and for $t=3,$ we get 
$F_{3}\left( \zeta \right) =v\left( \phi _{1},\phi _{2}\right) +v\left( \rho
_{1},\rho _{2}\right) e_{2}+v\left( \rho _{3},\rho _{4}\right) e_{4}.$ $%
F_{2}\left( \zeta \right) $ and $F_{3}\left( \zeta \right) $ are
hyperholomorphic, respectively octonionic hyperholomorphic function, from
Remark 2.1 and Theorem 3 from \cite{Xi-Zh-Li-05}.

For $t\geq 4$, using induction steps, supposing that $F_{t-1}\left( \zeta
\right) $ is a left $A_{t-1}$-hyperholomorphic function, we will prove that $%
F_{t}\left( \zeta \right) $ is $A_{t}$-hyperholomorphic. That means $%
D[F_{t}]=0.$ From relations (\ref{2.1}) and (\ref{2.2}), we have that 
\begin{equation*}
D[F_{t}]=\underset{k=0}{\overset{2^{t}-1}{\sum }}e_{k}\frac{\partial F_{t}}{%
\partial x_{k}}=\overset{2^{t-1}-1}{\underset{k=0}{\sum }}e_{k}\frac{%
\partial F_{t}}{\partial x_{k}}+\underset{k=2^{t-1}}{\overset{2^{t}-1}{\sum }%
}e_{k}\frac{\partial F_{t}}{\partial x_{k}}=
\end{equation*}%
\begin{equation*}
=D[F_{t-1}]+e_{2^{t-1}}\underset{k=0}{\overset{2^{t-1}-1}{\sum }}\overline{e}%
_{k}\frac{\partial F_{t}}{\partial x_{k+2^{t-1}}}.
\end{equation*}%
From induction steps, we obtain $D[F_{t-1}]=0.$ We will prove that $\underset%
{k=0}{\overset{2^{t-1}-1}{\sum }}\overline{e}_{k}\frac{%
\begin{array}{c}
\\ 
\partial F_{t}%
\end{array}%
}{%
\begin{array}{c}
\partial x_{2^{t-1}+k} \\ 
\end{array}%
}=0.$ This sum has $2^{t-1}$ terms. First two terms are: \newline
\begin{equation*}
\Big(\frac{\partial F_{t}}{\partial x_{2^{t-1}}}-e_{1}\frac{\partial F_{t}}{%
\partial x_{2^{t-1}+1}}\Big)=
\end{equation*}

\begin{equation*}
=\frac{\partial v}{\partial \rho _{2^{t-1}-1}}\frac{\partial \rho
_{2^{t-1}-1}}{\partial x_{2^{t-1}}}+\frac{\partial v}{\partial \rho
_{2^{t-1}}}\frac{\partial \rho _{2^{t-1}}}{\partial x_{2^{t-1}}}-e_{1}\left( 
\frac{\partial v}{\partial \rho _{2^{t-1}-1}}\frac{\partial \rho _{2^{t-1}-1}%
}{\partial x_{2^{t-1}+1}}+\frac{\partial v}{\partial \rho _{2^{t-1}}}\frac{%
\partial \rho _{2^{t-1}}}{\partial x_{2^{t-1}+1}}\right) =
\end{equation*}

\begin{equation*}
=\frac{\partial v}{\partial \rho _{2^{t-1}-1}}+\frac{\partial v}{\partial
\rho _{2^{t-1}}}\left( \frac{-1}{e_{1}}\right) -e_{1}\left( \frac{\partial v%
}{\partial \rho _{2^{t-1}-1}}\left( -e_{1}\right) +\frac{\partial v}{%
\partial \rho _{2^{t-1}}}\right) =
\end{equation*}

\begin{equation*}
=\frac{\partial v}{\partial \rho _{2^{t-1}-1}}+\frac{\partial v}{\partial
\rho _{2^{t-1}}}e_{1}-\frac{\partial v}{\partial \rho _{2^{t-1}-1}}-e_{1}%
\frac{\partial v}{\partial \rho _{2^{t-1}}}=0.
\end{equation*}%
Since $e_{1}^{2}=\gamma _{1},\,\,\gamma _{1}^{2}=1,$ $\frac{\partial v}{%
\partial \rho _{2^{t-1}-1}}$ and $\frac{\partial v}{\partial \rho _{2^{t-1}}}
$ can be written as $a_{2^{t-1}-1}\left( \zeta \right) +b_{2^{t-1}-1}\left(
\zeta \right) e_{1},$ respectively $a_{2^{t-1}}\left( \zeta \right)
+b_{2^{t-1}}\left( \zeta \right) e_{1}\ $ where $a_{2^{t-1}-1}\left( \zeta
\right) ,$ $b_{2^{t-1}-1}\left( \zeta \right) ,$ $a_{2^{t-1}}\left( \zeta
\right) $, $b_{2^{t-1}}\left( \zeta \right) $ are real valued functions.

\textit{Case 1}: $r<k.$ In the general case, we denote $%
T=2^{r}+2^{r+1}+...+2^{k}+2^{t-1}$ and $T_{1}=2^{r}+2^{r+1}+...+2^{k},$ for $%
r<k.$ We will compute the terms 
\begin{equation*}
-e_{T_{1}}\frac{\partial F_{t}}{\partial x_{T}}-e_{T_{1}+1}\frac{\partial
F_{t}}{\partial x_{T+1}}.
\end{equation*}%
We compute first $\frac{\partial F_{t}}{\partial x_{T}}.$ It results

\begin{equation*}
\frac{\partial F_{t}}{\partial x_{T}}=\biggr(...\left(\frac{\partial v}{%
\partial \rho _{T-1}} \frac{\partial \rho _{T-1}}{\partial x_{T}}+\frac{%
\partial v}{\partial \rho _{T}}\frac{\partial \rho _{T}}{\partial x_{T}}
\right)e_{2^{r}})e_{2^{r+1}})...e_{2^{k}})e_{2^{t-1}}=
\end{equation*}

\begin{equation*}
=\biggr(...\left(\frac{\partial v}{\partial \rho _{T-1}}+\frac{\partial v}{
\partial \rho _{T}}\frac{-1}{e_{1}}
\right)e_{2^{r}})e_{2^{r+1}})...e_{2^{k}})e_{2^{t-1}}=
\end{equation*}

\begin{equation*}
=\biggr(...\left( \frac{\partial v}{\partial \rho _{T-1}}+\frac{\partial v}{%
\partial \rho _{T}}e_{1}\right)
e_{2^{r}})e_{2^{r+1}})...e_{2^{k}})e_{2^{t-1}}.
\end{equation*}%
Since we can write $\frac{\partial v}{\partial \rho _{T-1}}$ under the form $%
a_{T-1}\left( \zeta \right) +b_{T-1}\left( \zeta \right) e_{1}$ and $\frac{%
\partial v}{\partial \rho _{T}}$ under the form $a_{T}\left( \zeta \right)
+b_{T}\left( \zeta \right) e_{1},$ where $a_{T-1}$, $b_{T-1}$, $a_{T}$, $%
b_{T}$ are real valued functions, using Proposition 2.2, we obtain:

\begin{equation*}
\frac{\partial F_{t}}{\partial x_{T}}=\biggr(...\left( \frac{\partial v}{%
\partial \rho _{T-1}}+\frac{\partial v}{\partial \rho _{T}}e_1 \right)
e_{2^{r}})e_{2^{r+1}})...e_{2^{k}})e_{2^{t-1}}=
\end{equation*}

\begin{equation*}
=(...(a_{T-1}(\zeta
)e_{2^{r}})e_{2^{r+1}})...e_{2^{k}})e_{2^{t-1}}+(...(b_{T-1}(\zeta
)e_{1})e_{2^{r}})e_{2^{r+1}})...e_{2^{k}})e_{2^{t-1}}+
\end{equation*}

\begin{equation*}
+(...(a_{T}(\zeta
)e_{1})e_{2^{r}})e_{2^{r+1}})...e_{2^{k}})e_{2^{t-1}}+(...(b_{T}(\zeta
)e_{1})e_{1})e_{2^{r}})e_{2^{r+1}})...e_{2^{k}})e_{2^{t-1}}=
\end{equation*}

\begin{equation*}
=a_{T-1}(\zeta )(-1)^{k-r+2}e_{T}+b_{T-1}(\zeta )(-1)^{k-r+3}e_{T+1}+
\end{equation*}

\begin{equation*}
+a_{T}(\zeta )(-1)^{k-r+3}e_{T+1}-b_{T}(\zeta )(-1)^{k-r+2}e_{T}.\medskip
\end{equation*}

Using Proposition 2.3, relation (\ref{2.6}), we compute $-e_{T_{1}}\frac{%
\partial F_{t}}{ \partial x_{T}}.$

\begin{equation*}
-e_{T_{1}}\frac{\partial F_{t}}{\partial x_{T}}=-e_{T_{1}}\biggr(%
a_{T-1}(\zeta )(-1)^{k-r+2}e_{T}+b_{T-1}(\zeta )(-1)^{k-r+3}e_{T+1}+
\end{equation*}

\begin{equation*}
+a_{T}(\zeta )(-1)^{k-r+3}e_{T+1}-b_{T}(\zeta )(-1)^{k-r+2}e_{T}\biggr)=
\end{equation*}

\begin{equation*}
=-\biggr(a_{T-1}(\zeta )(-1)^{k-r+2}(-1)^{k-r+1}e_{2^{i}}-b_{T-1}(\zeta
)(-1)^{k-r+3}(-1)^{k-r+1}e_{2^{i}+1}\biggr)-
\end{equation*}

\begin{equation*}
-\biggr(-a_{T}(\zeta )(-1)^{k-r+3}(-1)^{k-r+1}e_{2^{i}+1}-b_{T}(\zeta
)(-1)^{k-r+2}(-1)^{k-r+1}e_{2^{i}}\biggr)=
\end{equation*}

\begin{equation*}
=-\biggr(a_{T-1}(\zeta )(-1)^{2k-2r+3}e_{2^{i}}-b_{T-1}(\zeta
)(-1)^{2k-2r+4}e_{2^{i}+1}\biggr)-
\end{equation*}

\begin{equation*}
-\biggr(-a_{T}(\zeta )(-1)^{2k-2r+4}e_{2^{i}+1}-b_{T}(\zeta
)(-1)^{2k-2r+3}e_{2^{i}}\biggr).
\end{equation*}%
\medskip

Now, we compute $\frac{ \partial F_{t}}{ \partial x_{T+1}}.$ We obtain

\begin{equation*}
\frac{\partial F_{t}}{\partial x_{T+1}}=\biggr(...\biggr(\frac{\partial v}{%
\partial \rho _{T-1}}\frac{\partial \rho _{T-1}}{\partial x_{T+1}}+\frac{%
\partial v}{ \partial \rho _{T}}\frac{\partial \rho _{T}}{\partial x_{T+1}} %
\biggr)e_{2^{r}})e_{2^{r+1}})...e_{2^{k}})e_{2^{t-1}}=
\end{equation*}

\begin{equation*}
=\biggr(...\biggr(-\frac{\partial v}{\partial \rho _{T-1}}e_{1}+\frac{%
\partial v}{\partial \rho _{T}}\biggr )%
e_{2^{r}})e_{2^{r+1}})...e_{2^{k}})e_{2^{t-1}}.
\end{equation*}%
Since we can write $\frac{\partial v}{\partial \rho _{T-1}}$ under the form $%
a_{T-1}\left( \zeta \right) +b_{T-1}\left( \zeta \right) e_{1}$ and $\frac{%
\partial v}{\partial \rho _{T}}$ under the form $a_{T}\left( \zeta \right)
+b_{T}\left( \zeta \right) e_{1},$ where $a_{T-1}$, $b_{T-1}$, $a_{T}$, $%
b_{T}$ are real valued functions, using Proposition 2.2, we obtain:

\begin{equation*}
\frac{\partial F_{t}}{\partial x_{T+1}}=\biggr(...\biggr(-\frac{\partial v}{%
\partial \rho _{T-1}}e_1+\frac{\partial v}{\partial \rho _{T}} \biggr)%
e_{2^{r}})e_{2^{r+1}})...e_{2^{k}})e_{2^{t-1}}=
\end{equation*}

\begin{equation*}
=(...(-a_{T-1}(\zeta
)e_{1})e_{2^{r}})e_{2^{r+1}})...e_{2^{k}})e_{2^{t-1}}-(...(b_{T-1}(\zeta
)e_{1}e_{1})e_{2^{r}})e_{2^{r+1}})...e_{2^{k}})e_{2^{t-1}}+
\end{equation*}

\begin{equation*}
+(...(a_{T}(\zeta
))e_{2^{r}})e_{2^{r+1}})...e_{2^{k}})e_{2^{t-1}}+(...(b_{T}(\zeta
)e_{1}))e_{2^{r}})e_{2^{r+1}})...e_{2^{k}})e_{2^{t-1}}=
\end{equation*}

\begin{equation*}
=-a_{T-1}(\zeta )(-1)^{k-r+3}e_{T+1}+b_{T-1}(\zeta )(-1)^{k-r+2}e_{T}+
\end{equation*}

\begin{equation*}
+a_{T}(\zeta )(-1)^{k-r+2}e_{T}+b_{T}(\zeta )(-1)^{k-r+3}e_{T+1}.
\end{equation*}

\medskip Using Proposition 2.3, we compute $-e_{T_{1}+1}\frac{\partial F_{t}%
} {\partial x_{T+1}}.$

\begin{equation*}
-e_{T_{1}+1}\frac{\partial F_{t}}{\partial x_{T+1}}=-e_{T_{1}+1}\biggr(%
-a_{T-1}(\zeta )(-1)^{k-r+3}e_{T+1}+b_{T-1}(\zeta )(-1)^{k-r+2}e_{T}+
\end{equation*}

\begin{equation*}
+a_{T}(\zeta )(-1)^{k-r+2}e_{T}+b_{T}(\zeta )(-1)^{k-r+3}e_{T+1}\biggr)=
\end{equation*}

\begin{equation*}
=-\biggr(a_{T-1}(\zeta )(-1)^{k-r+3}(-1)^{^{k-r+1}}e_{2^{i}}-b_{T-1}(\zeta
)(-1)^{k-r+2}(-1)^{^{k-r+1}}e_{2^{i}+1}\biggr)-
\end{equation*}

\begin{equation*}
-\biggr(-a_{T}(\zeta )(-1)^{k-r+2}(-1)^{^{k-r+1}}e_{2^{i}+1}-b_{T}(\zeta
)(-1)^{k-r+3}(-1)^{^{k-r+1}}e_{2^{i}}\biggr)=
\end{equation*}

\begin{equation*}
=-\biggr(a_{T-1}(\zeta )(-1)^{2k-2r+4}e_{2^{i}}-b_{T-1}(\zeta
)(-1)^{2k-2r+3}e_{2^{i}+1}\biggr)-
\end{equation*}

\begin{equation*}
-\biggr(-a_{T}(\zeta )(-1)^{2k-2r+3}e_{2^{i}+1}-b_{T}(\zeta
)(-1)^{2k-2r+4}e_{2^{i}}\biggr).
\end{equation*}

Now, we can compute $-e_{T_{1}}\frac{ \partial F_{t}}{ \partial x_{T} }%
-e_{T_{1}+1}\frac{ \partial F_{t}}{ \partial x_{T+1}}.$ It results

\begin{equation*}
-e_{T_{1}}\frac{\partial F_{t}}{\partial x_{T}}-e_{T_{1}+1}\frac{\partial
F_{t}}{\partial x_{T+1}}=
\end{equation*}

\begin{equation*}
=-\biggr(a_{T-1}(\zeta )(-1)^{2k-2r+3}e_{2^{i}}-b_{T-1}(\zeta
)(-1)^{2k-2r+4}e_{2^{i}+1}\biggr)-
\end{equation*}

\begin{equation*}
-\biggr(-a_{T}(\zeta )(-1)^{2k-2r+4}e_{2^{i}+1}-b_{T}(\zeta
)(-1)^{2k-2r+3}e_{2^{i}}\biggr)-
\end{equation*}

\begin{equation*}
-\biggr(a_{T-1}(\zeta )(-1)^{2k-2r+4}e_{2^{i}}-b_{T-1}(\zeta
)(-1)^{2k-2r+3}e_{2^{i}+1}\biggr)-
\end{equation*}

\begin{equation*}
-\biggr(-a_{T}(\zeta )(-1)^{2k-2r+3}e_{2^{i}+1}-b_{T}(\zeta
)(-1)^{2k-2r+4}e_{2^{i}}\biggr)=0.
\end{equation*}%
\medskip

\textit{Case 2}: $r=k,$ we use Proposition 2.2 and Proposition 2.3, relation
(\ref{2.7}) and it easy to show that 
\begin{equation*}
-e_{2^{k}}\frac{\partial F_{t}}{\partial x_{T}}-e_{2^{k}+1}\frac{\partial
F_{t}}{\partial x_{T+1}}=0.
\end{equation*}%
$\Box \medskip $

\textbf{Remark\ 2.13. } The above proposition generalizes Theorem 3 from 
\cite{Xi-Zh-Li-05}. \medskip 
\begin{equation*}
\end{equation*}

\textbf{The Algorithm\medskip }

1) Input $t.$

2) Input functions $v,\phi _{1},\phi _{2}.$

3) For $i\in \{1,...,t-1\},$ $k\in \{1,...,i\},$ $r\in \{1,...,k-1\},$ \,
compute $M_{rki}=2^{r}+...+2^{k}+2^{i},$ $v\left( \rho _{M_{rki}-1},\rho
_{M_{rki}}\right) =\alpha _{M_{rki}}+\beta _{M_{rki}}e_{1}.$

4) For $i\in \{1,...,t-1\},$ $k\in \{1,...,i\},$ $r\in \{1,...,k-1\},$ 
\newline
\newline
-if $r<k,$ \thinspace we compute 
\begin{equation*}
(...\left( \alpha _{M_{rki}}+\beta _{M_{rki}}e_{1}\right)
e_{2^{r}})e_{2^{r+1}}...)e_{2^{k}})e_{2^{i}})=
\end{equation*}

\begin{equation*}
=\left( -1\right) ^{k-r+2}\left( \alpha _{M_{rki}}e_{M_{rki}}-\beta
_{M_{rki}}e_{M_{rki}-1}\right)
\end{equation*}%
\qquad\ \ \ \ \ \ \ \newline
-if $r=k\,,$ we compute 
\begin{equation*}
v\left( \rho _{2^{i}-1},\rho _{2^{i}}\right) e_{2^{i}}=(\alpha
_{2^{i}-1}+\beta _{2^{i}-1}e_{1})e_{2^{i}}=
\end{equation*}%
\begin{equation*}
=\alpha _{2^{i}-1}e_{2^{i}}+\beta _{2^{i}-1}e_{2^{i}+1}.
\end{equation*}

5) Output function 
\begin{equation*}
F_{t}\left( \zeta \right) =v\left( \phi _{1},\phi _{2}\right) +\underset{i=4}%
{\overset{t-1}{\sum }}\underset{k=1}{(\overset{i}{\sum }}\underset{r=1}{%
\overset{k-1}{(\sum }}\left( -1\right) ^{k-r+2}\left( \alpha
_{M_{rki}}\left( \zeta \right) e_{M_{rki}}-\beta _{M_{rki}}\left( \zeta
\right) e_{M_{rki}-1}\right) ))+
\end{equation*}%
\begin{equation*}
+\underset{i=1}{\overset{t-1}{\sum }}\left( \alpha _{2^{i}-1}\left( \zeta
\right) e_{2^{i}}+\beta _{2^{i}-1}\left( \zeta \right) e_{2^{i}+1}\right) .
\end{equation*}%
\medskip

\textbf{Conclusion.} In this paper, we generalized the notion of left $%
A_{t}- $holomorphic functions from quaternions to all algebras obtained by
the Cayley-Dickson process and we provided an algorithm to find examples of
left $A_{t}-$hyperholomorphic functions, using the \textit{shuffling}
procedure given by Bales in \cite{Ba-09}.

The theory of the right $A_{t}-$holomorphic functions and the theory of the
right $A_{t}-$hyperholomorphic functions are similarly to the corresponding
theories for the left functions and can be easy treated, using the above
ideas and procedures.

\begin{equation*}
\end{equation*}

This paper \ is supported by the grant of CNCS (Romanian National Council of
Research) PN-II-ID-WE-2012-4-169. 
\begin{equation*}
\end{equation*}

\bigskip

Cristina FLAUT

{\small Faculty of Mathematics and Computer Science,}

{\small Ovidius University,}

{\small Bd. Mamaia 124, 900527, CONSTANTA,}

{\small ROMANIA}

{\small http://cristinaflaut.wikispaces.com/}

{\small http://www.univ-ovidius.ro/math/}

{\small e-mail:}

{\small cflaut@univ-ovidius.ro}

{\small cristina\_flaut@yahoo.com}%
\begin{equation*}
\end{equation*}

Vitalii \ SHPAKIVSKYI

{\small Department of Complex Analysis and Potential Theory}

{\small \ Institute of Mathematics of the National Academy of Sciences of
Ukraine,}

{\small \ 3, Tereshchenkivs'ka st.}

{\small \ 01601 Kiev-4}

{\small \ UKRAINE}

{\small \ http://www.imath.kiev.ua/\symbol{126}complex/}

{\small \ e-mail: shpakivskyi@mail.ru}

\end{document}